\documentclass[10pt]{bmc_article}

\usepackage{cite} 
\usepackage{url}  
\usepackage{ifthen}  
\usepackage{multicol}   
\usepackage{amsmath, amsfonts, amssymb}
\newcommand{\R}{{{\mathbb R}}}

\def\includegraphics{}

\def\qed{\hbox to 0pt{}\hfill$\rlap{$\sqcap$}\sqcup$}

\newboolean{publ}

\newenvironment{bmcformat}{\baselineskip20pt\sloppy\setboolean{publ}{false}}{\baselineskip20pt\sloppy}

\begin{document}
\begin{bmcformat}

\title{On the Volterra property of a boundary problem with integral gluing condition for mixed parabolic-hyperbolic equation}

\author{Berdyshev A.S.$^1$%
         \email{Berdyshev A.S. - berdyshev@mail.ru},
       Cabada A.$^2$%
         \email{Cabada A. - alberto.cabada@usc.es},%
         Karimov E.T. \correspondingauthor$^3$%
         \email{Karimov E.T.\correspondingauthor - erkinjon@gmail.com}
         and
         Akhtaeva N.S. $^4$%
         \email{Akhtaeva N.S. - 260503@mail.ru}
      }

\address{%
    \iid(1)Kazakh National Pedagogical University named after Abai, Almaty, Kazakhstan\\
    \iid(2)University of Santiago de Compostela, Santiago de Compostela, Spain\\
    \iid(3)Institute of Mathematics, National University of Uzbekistan named after Mirzo Ulughbek, Tashkent, Uzbekistan\\
    \iid(4)Kazakh National Pedagogical University named after Abai, Almaty, Kazakhstan
}%

\maketitle

\begin{abstract}
In the present work we consider a boundary value problem with gluing conditions of integral form for parabolic-hyperbolic type equation. We prove that the considered problem has the Volterra property. The main tools used in the work are related to the method of the integral equations and functional analysis.
\end{abstract}

\ifthenelse{\boolean{publ}}{\begin{multicols}{2}}{}

\section*{Introduction}
	The theory of mixed type equations is one the principal parts of the general theory of partial differential equations. The interest for these kinds of equations arises intensively because of both theoretical and practical of their applications. Many mathematical models of applied problems require investigations of this type of equations.

	The actuality of the consideration of mixed type equations has been mentioned, for the first time,  by S. A. Chaplygin in 1902 in his famous work ``On gas streams'' \cite{chap}. The first fundamental results in this direction was obtained in 1920-1930 by F. Tricomi \cite{tric} and S. Gellerstedt \cite{gell}. The works of M. A. Lavrent'ev \cite{lavbit}, A. V. Bitsadze \cite{bit}, \cite{bit2}, F. I. Frankl \cite{fra}, M. Protter \cite{pro} and C. Morawetz \cite{mor}, have had a great impact in this theory, where outstanding theoretical results were obtained and pointed out important practical values of them. Bibliography of the main fundamental results on this direction can be found, among others, in the monographs of A. V.  Bitsadze \cite{bit2}, Y. M. Berezansky \cite{bere}, L. Bers \cite{bers}, M. S. Salakhitdinov and A. K. Urinov \cite{sal} and A. M. Nakhushev \cite{nak}.

	In most of the works devoted to the study of mixed type equations, the object of study was mixed elliptic-hyperbolic type equations. Comparatively, few results have been obtained on the study of mixed parabolic-hyperbolic type equations. However, this last type of equations  have also numerous applications in the real life processes (see \cite{ost} for an interesting example in mechanics). The reader can found a nice example given, for the first time, by Gelfand in \cite{gel}, and connected with the movement of the gas in a channel surrounded by a porous environment. Inside the channel the movement of gas was described by the wave equation and outside by the diffusion one. Mathematic models of this kind of problems arise in the study of electromagnetic fields, in heterogeneous environment, consisting of dielectric and conductive environment for modeling the movement of a little compressible fluid in a channel surrounded by a porous medium \cite{lei}. Here the wave equation describes the hydrodynamic pressure of the fluid in the channel, and the equation of filtration-pressure fluid in a porous medium. Similar problems arise in the study of the magnetic intensity of the electromagnetic field \cite{lei}.

	In the last few years, the investigations on local boundary value problems, for mixed equations in domains with non-characteristic boundary data, were intensively increased. We point out that the studies made on boundary value problems for equations of mixed type, in domains with deviation from the characteristics (with non-characteristic boundary), have been originated with the fundamental works of Bitsadze \cite{bit}, where the generalized Tricomi problem (Problem M) for an equation of mixed type is discussed.

	In the works \cite{salb} and \cite{sal}, the analog to the Tricomi problem for a modeled parabolic-hyperbolic equation, was investigated in a domain with non-characteristic boundary in a hyperbolic part. Moreover, the uniqueness of solution and the Volterra property of the formulated problem was proved.
	We also refer to the recent works devoted to the study of parabolic-hyperbolic equations \cite{asho}-\cite{khub}.

	In the last years, the interest for considering boundary value problems of parabolic-hyperbolic type, with integral gluing condition on the line of type changing, is increasing \cite{berk}, \cite{kapm}.

	 In the present work we study the analog to the generalized Tricomi problem with integral gluing condition on the line of type changing. We prove that the formulated problem has the Volterra property. The obtained result generalizes some previous ones from M. A. Sadybekov and G. D. Tajzhanova given in \cite{sadt}.

\section*{Formulation of the problem }
Let $\Omega  \subset {\R^2}$ be a domain, bounded at $y > 0$ by segments $A{A_0},{A_0}{B_0},B{B_0}$ of straight lines $x = 0,y = 1,x = 1$ respectively, and at  $y < 0$ by a monotone smooth curve $AC:\,y =  - \gamma (x),\,0 < x < l,\, 1/2 < l < 1,\,\gamma (0) = 0,\,l + \gamma (l) = 1$ and by the segment $BC:\,x - y = 1,\,l \leq x < 1$, which is the characteristic curve of the equation
$$
Lu = f(x,y),\eqno (1)
$$
where
$$
L\, u = \left\{ \begin{array}{l}
  {u_x} - {u_{yy}},\,\,\,\,y > 0, \hfill \\ \\
  {u_{xx}} - {u_{yy}},\,\,\,\,y < 0. \hfill \\
\end{array}  \right.\eqno (2)
$$

Now we state the problem that we will consider along the paper:

\textbf{Problem B.} To find a solution of the Eq.(1), satisfying boundary conditions
$$
{\left. {u(x,\,\,y)} \right|_{A{A_0} \cup {A_0}{B_0}}} = 0,\eqno (3)
$$
$$
{\left. {({u_x} - {u_y})} \right|_{AC}} = 0\eqno (4)
$$
and gluing conditions
$$
{u_x}(x, + 0) = {u_x}(x, - 0),\,\,\,\,\,{u_y}(x, + 0) = \alpha \,{u_y}(x,\, - 0) + \beta \int\limits_0^x {{u_y}(t,\, - 0)Q(x,t)\,dt} ,\,\,\,0 < x < 1,\eqno (5)
$$
where $Q$ is a given function such that $Q \in {C^1}\left( {\left[ {0,1} \right] \times \left[ {0,1} \right]} \right)$, and $\alpha ,\,\,\beta  \in \R$ satisfy ${\alpha ^2} + {\beta ^2} > 0$.

When the curve $AC$ coincides with the characteristic one $x + y = 0$, $\alpha  = 1$ and $\beta  = 0$, the problem B is just the Tricomi problem for parabolic-hyperbolic equation with non-characteristic line of type changing, which has been studied in \cite{ele}.

Regular solvability of the problem B with continuous gluing conditions  ($\alpha  = 1,\,\,\,\,\beta  = 0$) have been proved, for the first time, in \cite{ber}, and strong solvability of this problem was proved in the work \cite{sadt}.

Several properties, including the Volterra property of boundary problems for mixed parabolic-hyperbolic equations, have been studied in the works \cite{ber2}-\cite{berck}.

We denote the parabolic part of the mixed domain $\Omega$ as $\Omega_0$ and the hyperbolic part by $\Omega_1$.

A regular solution of the problem B in the domain $\Omega$ will be a function
$$
u \in C(\bar \Omega ) \cap {C^1}({\Omega _0} \cup AB)\cap {C^1}({\Omega _1} \cup AC \cup AB) \cap {C^{1,2}}({\Omega _0})\, \cap {C^{2,2}}({\Omega _1}),
$$
that satisfies Eq.(1) in the domains $\Omega_0$ and $\Omega_1$, the boundary conditions (3)-(4), and the gluing condition (5).

Regarding the curve $AC$, we assume that $x + \gamma (x)$ is monotonically increasing. Then, rewriting it by using the characteristic variables $\xi  = x + y$ and $\eta  = x - y$, we have that the equation of the curve $AC$ can be expressed as $\xi  = \lambda (\eta ),\,\,\,0 \leq \eta  \leq 1$.

\section*{Main result}

\textbf{Theorem 1.} Let $\gamma \in {C^1}[0,l]$ and $Q \in {C^1}\left( {\left[ {0,1} \right] \times \left[ {0,1} \right]} \right)$. Then for any function $f \in {C^1}(\bar \Omega )$, there exists a unique regular solution of the Problem B.

\verb"Proof:" By a regular solution of the problem B in the domain ${\Omega _1}$ we look for a function that fulfills the following expression
$$
u(\xi ,\,\,\eta ) = \frac{1}{2}\left[ {\tau (\xi ) + \tau (\eta ) - \int\limits_\xi ^\eta  {{\nu _1}(t)dt} } \right] - \int\limits_\xi ^\eta  {d{\xi _1}\int\limits_{{\xi _1}}^\eta  {{f_1}({\xi _1},\,{\eta _1})} d{\eta _1}},\eqno (6)
$$
where
$$
\xi  = x + y,\,\,\,\,\eta  = x - y,\,\,\,{f_1}(\xi ,\eta ) = \frac{1}{4}f(\frac{{\xi  + \eta }}{2},\,\frac{{\xi  - \eta }}{2}\,),\,\,\,\,\,\tau (x) = u(x,\, - 0),\,\,\,\,\,\,\,\,{\nu _1}(x) = {u_y}(x,\, - 0).\eqno (7)
$$

Based on (4) from (6), using the expressions on (7), we deduce that
$$
{\nu _1}(\eta ) = \tau '(\eta ) - 2\int\limits_{\lambda (\eta )}^\eta  {{f_1}({\xi _1},\,\eta )} d{\xi _1},\,\,\,0 \leq \eta  \leq 1.\eqno (8)
$$

By virtue of the unique solvability of the first boundary problem for the heat equation (1) satisfying condition (3), and the fact that $u\left( {x,0} \right) = \tau \left( x \right)$, its solution can be represented as
$$
u(x,y) = \int\limits_0^x {d{x_1}} \int\limits_0^1 {G(x - {x_1},\,y,{y_1})f({x_1},\,{y_1})} d{y_1} + \int\limits_0^x {{G_{{y_1}}}} (x - {x_1},y,0)\tau ({x_1})d{x_1},\eqno (9)
$$
where $\tau (0) = 0$ and $G\left( {x,y,{y_1}} \right)$ is the Green's function related to the first boundary problem, for the heat equation in a rectangle $A{A_0}{B_0}B$, which has the form \cite{tik}
$$
G(x,y,{y_1}) = \frac{1}{{2\sqrt {\pi x} }}\sum\limits_{n =  - \infty }^{ + \infty } {\left[ {\exp \left\{ { - \frac{{{{(y - {y_1} + 2n)}^2}}}{{4x}}} \right\} - \exp \left\{ { - \frac{{{{(y + {y_1} + 2n)}^2}}}{{4x}}} \right\}} \right]}.\eqno (10)
$$

Calculating the derivative  $\frac{{\partial u}}{{\partial y}}$  in  (9) and passing to the limit at $y \to 0$ we get
$$
{u_y}(x, + 0) =  - \int\limits_0^x {k(x - t){u_x}(t, + 0)dt}  + {F_0}(x),
$$
where				
$$
k(x) = \frac{1}{{\sqrt {\pi x} }}\sum\limits_{n =  - \infty }^{ + \infty } {{e^{ - \frac{{{n^2}}}{x}}}}  = \frac{1}{{\sqrt \pi  {x^{\frac{1}{2}}}}} + \widetilde k(x),\eqno (11)
$$
and
$$
{F_0}(x) = \int\limits_0^x {d{x_1}} \int\limits_0^1 {{{\left. {{G_y}(x - {x_1},y,{y_1})} \right|}_{y = 0}}f({x_1},{y_1})d{y_1}}.\eqno (12)
$$

Thus, the main functional relation between $\tau '(x)$ and ${\nu _0}(x) = {u_y}(x, + 0)$, reduced to the segment $AB$ from the parabolic part of the domain, imply that
$$
{\nu _0}(x) =  - \int\limits_0^x {k(x - t)\tau '(t)dt}  + {F_0}(x).\eqno (13)
$$

Suppose, in a first moment, that $\alpha  \ne 0$.  From (8) and (13), considering the gluing condition (5), we obtain the following integral equation regarding the function $\tau '(x)$:
$$
\tau '(x) + \int\limits_0^x {{k_1}(x,t)} \tau '(t)dt = {F_1}(x).\eqno (14)
$$

Here
$$
{k_1}(x,t) = \frac{1}{\alpha }\left[ {k(x - t) + \beta Q(x,t)} \right],\eqno (15)
$$
and
$$
{F_1}(x) = \frac{1}{\alpha }{F_0}(x) + 2\int\limits_{\lambda (x)}^x {{f_1}({\xi _1},x)d{\xi _1}}  + \frac{{2\beta }}{\alpha }\int\limits_0^x {Q(x,t)dt\int\limits_{\lambda (t)}^t {f({\xi _1},t)d{\xi _1}} }. \eqno (16)
$$

Hence, the problem B  is equivalent, in the sense of unique solvability, to the second kind Volterra integral equation (14).

The restrictions imposed on the functions $\gamma$, $Q$, and the right hand of the Eq.(1) guarantee that, by virtue of (11) and (15), the kernel ${k_1}(x,t)$ is a kernel with weak singularity. So, we have that Eq.(14) has a unique solution and $\tau' \in {C^1}(0,1)$. Since $\tau(0)=0$, we deduce the uniqueness of the function $\tau$. Eq.(8) gives us the uniqueness of function $\nu_1$ and, as consequence, we deduce, from Eq.(6), the uniqueness of solution of problem B when $\alpha  \ne 0$.

Consider now the other case, i.e. $\alpha  = 0$  and $\beta  \ne 0$.

From functional relations (8) and (13), and taking gluing condition (5) into account at $\alpha  = 0$, we have
$$
 - \int\limits_0^x {k(x - t)\tau '} (t)dt + {F_0}(x) = \beta \int\limits_0^x {[\tau '(t) - 2\int\limits_{\lambda (t)}^t {{f_1}({\xi _1},t)d{\xi _1}} ]Q(x,t)dt}
 $$
or, which is the same,
$$
\int\limits_0^x {\tau '(t)[k(x - t) + \beta Q(x,t)]dt}  = {F_0}(x) + 2\beta \int\limits_0^x {dt\int\limits_{\lambda (t)}^t {Q(x,t){f_1}({\xi _1},t)d{\xi _1}} .}
$$
Considering the representation of $k(x - t)$, the previous equation can be rewritten as follows
$$
\int\limits_0^x {\frac{\tau '(t)dt}{{\left(x - t\right)}^{1/2}}}  = \sqrt \pi  \left[ F_0(x) + 2\beta \int\limits_0^x {dt} \int\limits_{\lambda (t)}^t {Q(x,t){f_1}(\xi _1,t)d{\xi _1}}  - \int\limits_0^x {\tau'(t)\left(\widetilde{k}(x - t)+ \beta Q(x,t)\right) dt} \right].\eqno (17)
$$
Since Eq.(17) is the Abel's equation, it can be solved and so we arrive at the following identity:
\begin{eqnarray*}
 \tau '(x) &=& \displaystyle \frac{{{F_0}(0)}}{{\sqrt {\pi x} }} + \frac{1}{{\sqrt \pi  }}\left\{ {\int\limits_0^x {\frac{{{F'_0}(t)dt}}{{\sqrt {x - t} }} + 2\beta \int\limits_0^x {\frac{{dt}}{{\sqrt {x - t} }}} } \frac{\partial }{{\partial t}}\int\limits_0^t {dz\int\limits_{\lambda (z)}^z {Q(t,z){f_1}({\xi _1},z)d{\xi _1}} } } \right. \hfill \\
  && - \left. {\int\limits_0^x {\frac{{dt}}{{\sqrt {x - t} }}} \frac{\partial }{{\partial t}}\int\limits_0^t {\tau (z)\left[ {\widetilde k(x - t) + \beta Q(t,z)} \right]dz} } \right\}. \hfill \\
\end{eqnarray*}

Considering ${F_0}(0) = 0$, after some simplifications we get
$$
\tau '(x) + \int\limits_0^x {{K_0}(x,z)} \tau '(z)dz = {F_2}(x),\eqno (18)
$$
where
$$
{K_0}(x,z) = \frac{1}{{\sqrt \pi  }}\left\{ {\frac{{Q(z,z)}}{{\sqrt {x - z} }} + \int\limits_0^{x - z} {{{(x - t)}^{ - \frac{1}{2}}}} \frac{\partial }{{\partial t}}\left[ {\widetilde k(t - z) + \beta Q(t,z)} \right]dt} \right\},
$$

\begin{eqnarray*}
\hspace{3cm}  {F_2}(x) &=&\frac{1}{{\sqrt \pi  }}\int\limits_0^x {d{x_1}} \int\limits_0^1 {\left[ {\int\limits_0^{x - {x_1}} {\frac{{{G_{yt}}(t,{y_1},0)}}{{\sqrt {x - {x_1} - t} }}dt} } \right]f({x_1},{y_1})d{y_1}} \hfill \\
    & & \hspace{11cm} (19)\\
&&+  \frac{{2\beta }}{{\sqrt \pi  }}\int\limits_0^x {d{\eta _1}} \int\limits_{\lambda ({\eta _1})}^{{\eta _1}} {\left[ {\frac{{Q({\eta _1},{\eta _1})}}{{\sqrt {x - {\eta _1}} }} + \int\limits_0^{x - {\eta _1}} {\frac{{{Q_t}(t,{\eta _1})}}{{\sqrt {x - {\eta _1} - t} }}dt} } \right]} {f_1}({\xi _1},{\eta _1})d{\xi _1}. \hfill \\
\end{eqnarray*}

Since the kernel ${K_0}(x,z)$ has a weak singularity, then Eq.(18) has a unique solution, and it can be represented as
$$
\tau '(x) = {F_2}(x) + \int\limits_0^x {R(x,z)} {F_2}(z)dz,\eqno (20)
$$
where $R(x,z)$ is the resolvent kernel of (18).

As consequence, arguing as in the case $\alpha\neq 0$, we deduce, from Eq.(6), the uniqueness of solution of Problem B for $\alpha=0$ and $\beta\neq 0$, and the result is proved. \qed

In the sequel, we will deduce the exact expression of the integral kernel related to the unique solution of Problem B.

To this end, we suppose, at the beginning, that $\alpha\neq 0$. Note that the unique solution of Eq.(14) can be represented as
$$
\tau '(x) = \int\limits_0^x {\Gamma (x,t)} {F_1}(t)dt + {F_1}(x),\eqno (21)
$$
where $\Gamma (x,t)$ is the resolvent kernel of the Eq.(14), and it is given by the recurrence formula:
$$
\Gamma (x,t) = \sum\limits_{n = 1}^\infty  {{{( - 1)}^n}{k_{1n}}(x,t),\,\,\,\,{k_{11}}(x,t)}  = {k_1}(x,t),\,\,\,{k_{1n + 1}}(x,t) = \int\limits_0^x {{k_1}(x,z){k_{1n}}(t,z)} dz.
$$

From (21), taking $\tau (0) = 0$ into account, we have that
$$
\tau (x) = \int\limits_0^x {{\Gamma _1}(x,t)} {F_1}(t)dt,
$$
where  ${\Gamma _1}(x,t) = 1 + \int\limits_t^x {\Gamma (z,t)} dz$.

From the formula (6), and considering (8), one can easily deduce that
$$
u(\xi,\eta) = \tau (\xi ) + \int\limits_\xi ^\eta  {d{\eta _1}} \int\limits_{\lambda ({\eta _1})}^{\eta_1}  {f({\xi _1},{\eta _1})d{\xi _1}} .
\eqno (22)
$$

Substituting the representation of $\tau (x)$ into (22) and considering (12) and (16), after some evaluations we get
\begin{eqnarray*}
 \hspace{2cm} u(x,y) &=& \frac{1}{\alpha }\int\limits_0^\xi  {d{x_1}} \int\limits_0^1 {{G_1}(\xi  - {x_1},{y_1})f({x_1},{y_1})d{y_1}}  + 2\int\limits_0^\xi  {d{\eta _1}} \int\limits_{\lambda ({\eta _1})}^{{\eta _1}} {{\Gamma _1}(\xi ,{\eta _1}){f_1}({\xi _1},{\eta _1})d{\xi _1}}  \hfill \\
  &&   \hspace{12cm} (23)\\
  && + \frac{{2\beta }}{\alpha }\int\limits_0^\xi  {d{\eta _1}} \int\limits_{\lambda ({\eta _1})}^{{\eta _1}} {{G_0}(\xi  - {\eta _1},{\eta _1}){f_1}({\xi _1},{\eta _1})d{\xi _1}}  + \int\limits_\xi ^\eta  {d{\eta _1}} \int\limits_{\lambda ({\eta _1})}^{\eta_1}  {f({\xi _1},{\eta _1})d{\xi _1},}
\end{eqnarray*}
where
$$
{G_1}(x,{y_1}) = \int\limits_0^x {{\Gamma _1}(x,t){G_y}(t,{y_1},0)dt},
$$
and
$$
{G_0}(x,\eta ) = \int\limits_0^x {Q(z + \eta ,\eta ){\Gamma _1}(x,z)dz}.
$$

In an analogous way,  substituting the representation of $\tau (x)$ into (9), we have
\begin{eqnarray*}
 \hspace{2cm}
  u(x,y) &=& \int\limits_0^x {d{x_1}} \int\limits_0^1 {{G_2}(x - {x_1},y,{y_1})f({x_1},{y_1})d{y_1}}  + 2\int\limits_0^x {d{\eta _1}} \int\limits_{\lambda ({\eta _1})}^{{\eta _1}} {{G_1}(x - {\eta _1},y){f_1}({\xi _1},{\eta _1})d{\xi _1}}  \hfill \\
   &&   \hspace{12cm} (24)\\
  &&   + \frac{{2\beta }}{\alpha }\int\limits_0^x {d{\eta _1}} \int\limits_{\lambda ({\eta _1})}^{{\eta _1}} {{G_{01}}(x - {\eta _1},{\eta _1}){f_1}({\xi _1},{\eta _1})d{\xi _1}} , \hfill \\
\end{eqnarray*}
where
$$
{G_2}(x,y,{y_1}) = G(x,y,{y_1}) + \frac{1}{\alpha }\int\limits_0^x {{G_1}(t,{y_1}){G_y}(x - t,y,0)dt}
$$
and
$$
{G_{01}}(x,{\eta _1}) = \int\limits_0^x {{G_y}({x_1},{y_1},0){G_0}(x - {x_1},{\eta _1})d{x_1}} .
$$

From (23) and (24), we arrive at the following expression
$$
u(x,y) = \iint\limits_\Omega  {{K_{\alpha \beta }}(x,y,{x_1},{y_1})f({x_1},{y_1})}d{x_1}d{y_1},
$$			
where
\begin{eqnarray*}
  {K_{\alpha \beta }}(x,y,{x_1},{y_1}) &=& \theta (y)\left\{ {\frac{{}}{{}}\theta ({y_1})\theta (x - {x_1}){G_2}(x - {x_1},y,{y_1}) + \theta ( - {y_1})\theta (x - {\eta _1})\left[ {\frac{{}}{{}}{G_1}(x - {\eta _1},y) } \right.} \right. \hfill \\
   && \left. +{\left. {\frac{{2\beta }}{\alpha }{G_{01}}(x - {\eta _1},{\eta _1})} \right]} \right\} + \theta ( - y)\left\{ {\theta ({y_1})\theta (\xi  - {x_1}){G_1}(\xi  - {x_1},{y_1})} \right.  \hfill \\
 &&
+\left. \theta ( - {y_1})\left[ {\frac{1}{2}\theta (\eta  - {\eta _1})\theta ({\eta _1} - \xi )\theta (\xi  - {\xi _1})  } \right.    {+\left. {\theta (\xi  - {\eta _1})\left[ {{\Gamma _1}(\xi ,{\eta _1}) + \frac{\beta }{\alpha }{G_0}(\xi  - {\eta _1},{\eta _1})} \right]} \right]} \right\}. \hfill
\end{eqnarray*}

Here
$$
\theta (y)  = \left\{ \begin{array}{l}
  1,\,\,\,\,y > 0, \hfill \\ \\
  0,\,\,\,\,y < 0. \hfill \\
\end{array}  \right.
$$

When $\alpha=0$ and $\beta\neq 0$, by using a similar algorithm, we conclude that
$$
u(x,y) = \iint\limits_\Omega  {{K_{0\beta }}(x,y,{x_1},{y_1})f({x_1},{y_1})}d{x_1}d{y_1},
$$  				
where
\begin{eqnarray*}
  {K_{0\beta }}(x,y,{x_1},{y_1}) &=& \theta (y)\left\{ {\theta ({y_1})\theta (x - {x_1})\left[ {G(x - {x_1},y,{y_1}) + {G_4}(x - {x_1},{x_1},y,{y_1})} \right]} \right. \hfill \\
&&  \left. { + \theta ( - {y_1})\theta (x - {\eta _1}){G_5}(x - {\eta _1},y,{\eta _1})} \right\} + \theta ( - y)\left\{ {\frac{{}}{{}}\theta ({y_1})\theta (\xi  - {x_1}){G_3}(\xi ,{x_1},{y_1})} \right. \hfill \\
 && \left. { + \theta ( - {y_1})\left[ {\theta (\xi  - {\eta _1}){Q_1}(\xi ,\,{\eta _1}) + \frac{1}{2}\theta (\eta  - {\eta _1})\theta ({\eta _1} - \xi )\theta (\xi  - {\xi _1})} \right]} \right\}, \hfill
\end{eqnarray*}
with
\begin{eqnarray*}
  {G_3}(x,{x_1},{y_1}) &=& \frac{1}{{\sqrt \pi  }}\int\limits_0^{x - {x_1}} {\left[ {\frac{{{G_y}(z,{y_1},0)}}{{\sqrt z }} + \int\limits_0^z {\left\{ {\frac{{{G_{ys}}(s,{y_1},0)}}{{\sqrt {z - s} }} + R(z + {x_1},s + {x_1})\left[ {\frac{{{G_y}(s,{y_1},0)}}{{\sqrt s }}} \right.} \right.} } \right.}  \hfill \\
 && \left. {\left. {\left. { + \int\limits_0^s {\frac{{{G_{yt}}(t,{y_1},0)}}{{\sqrt {s - t} }}dt} } \right]} \right\}ds} \right]dz, \hfill
\end{eqnarray*}

\begin{eqnarray*}
  {Q_1}(x,{\eta _1}) &=& \frac{{2\beta }}{{\sqrt \pi  }}\int\limits_0^{x - {\eta _1}} {\left[ {\frac{{Q({\eta _1},{\eta _1})}}{{\sqrt z }} + \int\limits_0^z {\left\{ {\frac{{{Q_s}(s,{\eta _1})}}{{\sqrt {z - s} }} + R(z + {\eta _1},s + {\eta _1})\left[ {\frac{{Q({\eta _1},{\eta _1})}}{{\sqrt s }}} \right.} \right.} } \right.}  \hfill \\
 && \left. {\left. {\left. { + \int\limits_0^s {\frac{{{Q_t}(t,{\eta _1})}}{{\sqrt {s - t} }}dt} } \right]} \right\}ds} \right]dz, \hfill \\
\end{eqnarray*}
$$
{G_4}(t,{x_1},{y_1},y) = \int\limits_0^t {{G_3}(s,{x_1},{y_1}){G_y}(t - s,y,0)ds}
$$
and
$$
{G_5}(x,y,{\eta _1}) = \int\limits_0^x {{G_y}(x - {x_1},y,0){Q_1}({x_1},{\eta _1})d{x_1}} .
$$

Thus we have partially proved the following lemma.

\textbf{Lemma 1.} The unique regular solution of Problem B can be represented as follows
$$
u\left( {x,y} \right) = \iint\limits_\Omega  {K\left( {x,y,{x_1},{y_1}} \right)f\left( {{x_1},{y_1}} \right)d{x_1}d{y_1}},\,\,(x,y)\in \Omega,\eqno (25)
$$
where $K\left( {x,y,{x_1},{y_1}} \right) \in {L_2}\left( {\Omega  \times \Omega } \right)$ and
$$K\left( {x,y,{x_1},{y_1}} \right) = {K_{\alpha \beta }}\left( {x,y,{x_1},{y_1}} \right), \quad \mbox{ if $\alpha  \ne 0,$}$$
$$K\left( {x,y,{x_1},{y_1}} \right) = {K_{0\beta }}\left( {x,y,{x_1},{y_1}} \right),  \quad \mbox{ if $\alpha  = 0.$}$$

\verb"Proof": Expression (25) has been proved before. Let's see that $K\left( {x,y;{x_1},{y_1}} \right) \in {L_2}\left( {\Omega  \times \Omega } \right)$.

Note that in the kernel defined in (25), all the items are bounded except the first one. So, we only need to prove that
$$
\theta \left( y \right)\theta \left( {{y_1}} \right)\theta \left( {x - {x_1}} \right)G\left( {x - {x_1},y,{y_1}} \right) \in {L_2}\left( {\Omega  \times \Omega } \right).
$$

From the representation of the Green's function $G\left( {x - {x_1},y,{y_1}} \right)$ given in  (10), it follows that, for the aforementioned aim, it is enough to prove that (for $n = 0$):
$$
B\left( {x - {x_1},y,{y_1}} \right) = \theta \left( y \right)\theta \left( {{y_1}} \right)\theta \left( {x - {x_1}} \right)\frac{1}{{2\sqrt {\pi \left( {x - {x_1}} \right)} }}
\, \left[ {\exp \left\{ { - \frac{{{{\left( {y - {y_1}} \right)}^2}}}{{4\left( {x - {x_1}} \right)}}} \right\} - \exp \left\{ { - \frac{{{{\left( {y + {y_1}} \right)}^2}}}{{4\left( {x - {x_1}} \right)}}} \right\}} \right]
$$
is bounded.

First, note that
$$
   B\left( x - x_1,y,y_1 \right) \leq \frac{1}{2\sqrt {\pi \left( x - x_1 \right)} }e^{ - \frac{\left( y - y_1\right)^2}{4\left( x - x_1\right)}}.
$$

Using this fact, we deduce that
\begin{eqnarray*}
\left\| B \right\|_{{L_2}\left( {\Omega  \times \Omega } \right)}^2 &=& \int\limits_0^1 {dx\int\limits_0^1 {dy\int\limits_0^x {d{x_1}\int\limits_0^1 {{{\left| {B\left( {x - {x_1},y,{y_1}} \right)} \right|}^2}d{y_1} } } } }\\
& = & \int\limits_0^1 {dy} \int\limits_0^1 d{y_1}\int\limits_0^1 dx\int\limits_0^x {\left| B\left( x,y,{y_1} \right) \right|^2 d{x_1} } \\
&\leq & \int\limits_0^1 {dy} {\int\limits_0^1 {d{y_1}\int\limits_0^1 {\left| {B\left( {x,y,{y_1}} \right)} \right|^2} } }dx \leq \frac{1}{{4\pi }}\int\limits_0^1 {dy\int\limits_0^1 {d{y_1}} \int\limits_0^1 {\frac{1}{x}{e^{ - \frac{{{{\left( {y - {y_1}} \right)}^2}}}{{4x}}}}dx } }
\\ &=&
 \frac{1}{4\pi }\int\limits_0^1{dy}\int\limits_0^1 {\frac{dx}{x}}\int\limits_0^1 {e^{ - \frac{\left(y -y_1\right)^2}{4x}}dy_1.}
\end{eqnarray*}

By means of the change of variables $\frac{y - y_1}{2\sqrt x } = y_2$, we get that this last expression is less than or equals to the following one
$$
\frac{1}{{4\pi }}\int\limits_0^1 {dy} \int\limits_0^1 {\frac{{dx}}{x}} \int\limits_{\frac{{y - 1}}{{2\sqrt x }}}^{\frac{y}{{2\sqrt x }}} {{e^{ - y_2^2}}2\sqrt x d{y_2}}
\leq \frac{1}{2\pi}\int\limits_0^1 {dy}\int\limits_0^1 {\frac{dx}{\sqrt {x} }}\int\limits_{-\infty}^{+\infty} {e^{-y_2}dy_2} = \frac{1}{\sqrt \pi}.
$$

As consequence, $K\left( {x,y;{x_1},{y_1}} \right) \in {L_2}\left( {\Omega  \times \Omega } \right)$ and Lemma 1 is completely proved. \qed

Define now
$$
F_{\alpha\beta}(x)=\left\{
\begin{array}{l}
F_1(x),\,\,\alpha\neq 0\\
F_2(x),\,\,\alpha=0.\\
\end{array}
\right.
$$
We have the following regularity result for this function:

\textbf{Lemma 2.} If $f \in {C^1}(\overline{ \Omega} ),\,\,f(0,0) = 0$ and $Q \in {C^1}\left( {\left[ {0,1} \right] \times \left[ {0,1} \right]} \right)$, then $F_{\alpha \beta } \in {C^1}[0,1]$ and $F_{\alpha \beta }\left( 0 \right) = 0$.

\verb"Proof:" Using the explicit form of the Green's function given in (10), it is not complicated to prove that function ${F_{\alpha \beta }}$, defined by formulas (16) and (19), belongs to the class of functions ${C^1}[0,1]$ and $F_{\alpha \beta }\left( 0 \right) = 0$. \\
Lemma 2 is proved. \qed

\textbf{Lemma 3.}  Suppose that  $Q \in {C^1}\left( {\left[ {0,1} \right] \times \left[ {0,1} \right]} \right)$ and $f \in {L_2}\left( \Omega  \right)$, then ${F_{\alpha \beta }} \in {L_2}(\Omega )$ and
$$
{\left\| {{F_{\alpha \beta }}} \right\|_{{L_2}\left( {0,1} \right)}} \leq C{\left\| f \right\|_{L_2(\Omega)}}.\eqno (26)
$$

\verb"Proof:"  Consider the following problem in ${\Omega _0}$:
$$
{\omega _x} - {\omega _{yy}} = f(x,y),\,\,\,\,\,\,\,\,\,\,\,\,\,\,{\left. \omega  \right|_{A{A_0} \cup {A_0}{B_0} \cup AB}} = 0.\eqno (27)
$$

Obviously, we have that ${F_0}(x) = \mathop {\lim }\limits_{y \to 0} {\omega _y}(x,y)$.

First, note that it is known \cite{mix} that problem (27) has a unique solution $\omega\in W_2^{1,2}({\Omega _0})$, and it satisfies the following inequality
$$
\left\| \omega  \right\|_{{L_2}({\Omega _0})}^2 + \left\| {{\omega _x}} \right\|_{{L_2}({\Omega _0})}^2 + \left\| {{\omega _y}} \right\|_{{L_2}({\Omega _0})}^2 + \left\| {{\omega _{yy}}} \right\|_{{L_2}({\Omega _0})}^2 \leq C\left\| f \right\|_{{L_2}({\Omega _0})}^2. \eqno (28)
$$

Using now the obvious equality
$$
{\omega _y}(x,0) = {\omega _y}(x,y) - \int\limits_0^y {{\omega _{yy}}(x,t)dt},
$$
we have that
$$
\left\| {{\omega _y}(\cdot,0)} \right\|_{{L_2}(0,1)}^2 = {\int\limits_0^1 {\left| {{\omega _y}(x,0)} \right|^2} }dx = \int\limits_0^1 {dy{{\int\limits_0^1 {\left| {{\omega _y}(x,0)} \right|^2} }}dx}  \leq C\left[ {\left\| {{\omega _y}} \right\|_{{L_2}({\Omega _0})}^2 + \left\| {{\omega _{yy}}} \right\|_{{L_2}({\Omega _0})}^2} \right].\eqno (29)
$$

From (28) and (29) we obtain
$$
{\left\| {{F_0}} \right\|_{{L_2}(0,1)}} = {\left\| {{\omega _y}\left( {\cdot,0} \right)} \right\|_{{L_2}(0,1)}} \leq C{\left\| f \right\|_{{L_2}({\Omega _0})}}.\eqno (30)
$$

Now, by virtue of the conditions of Lemma 3 and the representations (16) and (19), from expression (30) and the Cauchy-Bunjakovskii inequalities, we get the estimate (26) and conclude the proof. \qed

Denote now $\| \cdot \|_l$ as the norm of the Sobolev space $H^l(\Omega)\equiv W_2^l(\Omega)$ with $W_2^0(\Omega)\equiv L_2(\Omega)$.

\textbf{Lemma 4.} Let $u$ be the unique regular solution of Problem B. Then the following estimate holds:
$$
{\left\| u \right\|_1} \leq c{\left\| f \right\|_0}.\eqno (31)
$$
Here $c$ is a positive constant that does not depend on $u$.

\verb"Proof:" By virtue of Lemma 3, and from (20) and (21), we deduce that
$$
{\left\| {\tau '} \right\|_{{L_2}\left( {0,1} \right)}} \leq C{\left\| {{F_{\alpha \beta }}} \right\|_{{L_2}\left( {0,1} \right)}} \leq C{\left\| f \right\|_0}.
$$
The result follows from expression (22). \qed

\textbf{Definition 1}. We define the set $W$ as the set of all the regular solutions of Problem B.

A function $u \in {L_2}\left( \Omega  \right)$ is said to be a strong solution of Problem B, if there exists a functional sequence $\left\{ {{u_n}} \right\} \subset W$, such that ${u_n}$ and $L{u_n}$ converge in ${L_2}\left( \Omega  \right)$ to $u$ and $f$ respectively.

Define $\mathbb{L}$ as the closure of the differential operator $\mathbb{L}:\,W\rightarrow {L_2}\left( \Omega  \right)$, given by expression (2).

Note that, according to the definition of the strong solution, the function $u$ will be a strong solution of Problem B if and only if $u \in D\left( \mathbb{L} \right)$.

Now we are in a position to prove the following uniqueness result for strong solutions.

\textbf{Theorem 2.} For any function $Q \in {C^1}\left( {\left[ {0,1} \right] \times \left[ {0,1} \right]} \right)$ and $f \in {L_2}\left( \Omega  \right)$, there exists a unique strong solution $u$ of Problem B. Moreover $u\in W_2^1\left( \Omega  \right) \cap W_{x,y}^{1,2}\left( {{\Omega _1}} \right) \cap C\left( {\overline \Omega } \right)$, satisfies inequality (31) and it is given by the expression (25).

\verb"Proof:" Let $C_{0}^{1}\left( {\overline{\Omega }} \right)$ be the set of the $C^1\left(\overline{\Omega}\right)$ functions that vanish in a neighborhood of $\partial \Omega$  ($\partial \Omega$ is a boundary of the domain $\Omega$). Since $C_{0}^{1}\left( {\overline{\Omega }} \right)$ is dense in ${L_2}\left( \Omega  \right)$, we have that for any function $f \in {L_2}\left( \Omega  \right)$,  there exist a functional sequence ${{f}_{n}}\in C_{0}^{1}\left( {\overline{\Omega }} \right)$, such that $\left\| {{f_n} - f} \right\| \to 0$,  as $n \to \infty$.

It is not difficult to verify that if $f_n\in C_0^1\left( {\overline{\Omega }} \right)$ then $F_{\alpha \beta n}\in C^1(\left[ 0,1 \right])$ (with obvious notation). Therefore equations (14) and (18)  can be considered as a second kind Volterra integral equations in the space $C^1(\left[ 0,1 \right])$. Consequently, we have that ${\tau '_n}\left( x \right) = {u_{nx}}\left( {x,0} \right) \in {C^1}\left[ {0;1} \right]$. Due to the properties of the solutions of the boundary value problem for the heat equation in $\Omega_0$ and the Darboux problem, by using the representations (6) and (9), we conclude that ${u_n} \in W$ for all $f_n\in C_{0}^{1}\left( {\overline{\Omega }} \right)$.

By virtue of the inequality (31) we get
$$
{\left\| {{u_n} - u} \right\|_1} \leq c{\left\| {{f_n} - f} \right\|_0} \to 0.
$$
Consequently,  $\left\{ {{u_n}} \right\}$ is a sequence of strong solutions, hence, Problem B is strongly solvable for all right hand $f\in L_2(\Omega)$, and the strong solution belongs to the space $W_2^1\left( \Omega  \right) \cap W_{x,y}^{1,2}\left( {{\Omega _1}} \right) \cap C\left( {\overline \Omega } \right)$. Thus, Theorem 2 is proved. \qed

Consider now, for all $n=2,3, \ldots$, the sequence of kernels given by the recurrence formula
$$
{K_n}(x,y;\,{x_1},{y_1}) = \iint\limits_\Omega  {K(x,y;\,{x_2},{y_2}){K_{(n - 1)}}({x_2},{y_2},{x_1},{y_1})d{x_2}d{y_2}\,\,},
$$
with
$$
{K_1}(x,y;\,{x_1},{y_1}) = K(x,y;\,{x_1},{y_1}),
$$
and $K$ defined in Lemma 1.

\textbf{Lemma 5. } For the iterated kernels ${K_n}(x,y;\,{x_1},{y_1})$ we have the following estimate:
$$
\left| {{K_n}(x,y;\,{x_1},{y_1})} \right| \leq {(\sqrt {\pi \,} M)^n}{\left( {\frac{3}{2}} \right)^{n - 1}}\,\frac{{{{(x - {x_1})}^{\frac{n}{2} - 1}}}}{{\Gamma \left( {\frac{n}{2}} \right)}}\,\,,\,\,\,\,\,n = 1,2,3...,\eqno (32)
$$
where $
M=\underset{\begin{smallmatrix}
 \left( x,y \right)\in \Omega  \\
 \left( {{x}_{1}},{{y}_{1}} \right)\in \Omega
\end{smallmatrix}}{\mathop{\max }}\,\left| \sqrt{x-{{x}_{1}}}\,K\left( x,y;\,{{x}_{1}},{{y}_{1}} \right) \right|
$ and $\Gamma$ is the Gamma-function of Euler.

\verb"Proof:" The proof will be done by induction in $n$.

Taking the representation of the Green's function given in (10) into account, and from the representation of the kernel $K(x,y;\,{x_1},{y_1})$ at $n = 1$, the inequality (32)
$$
\left| {{K_1}(x,y;\,{x_1},{y_1})} \right| \leq M{(x - {x_1})^{ - \frac{1}{2}}}
$$
is automatically deduced.

Let (32) be valid for $n = k - 1$. We will prove the validity of this formula for $n = k$. To this end, by using inequality (32), at $n = 1$ and $n = k - 1$, we have that

\begin{eqnarray*}
  \left| K_k(x,y;x_1,y_1) \right| &=& \left| \iint\limits_\Omega  K(x,y;x_2,y_2)K_{(k - 1)}(x_2,y_2,x_1,y_1)dx_2dy_2 \right|  \hfill \\
 & \leq &   \iint\limits_\Omega  {\left| K(x,y;x_2,y_2) \right| \, \left| K_{(k - 1)}(x_2,y_2;x_1,y_1) \right|dx_2dy_2}  \hfill \\
 & \leq &   \iint\limits_\Omega  {\theta (x - x_2)M(x - x_2)^{ - \frac{1}{2}}} \, \theta (x_2 - x_1) \, (\sqrt \pi   {M})^{k - 1}\left(\frac{3}{2}\right)^{k - 2} \frac{(x_2 - x_1)^{\frac{k}{2} - \frac{3}{2}}}{\Gamma (\frac{k - 1}{2})}dx_2dy_2  \hfill \\
 & \leq & {M^k}{(\sqrt \pi  )^{k - 1}}{\left(\frac{3}{2}\right)^{k - 2}} \, \frac{1}{{\Gamma (\frac{{k - 1}}{2})}}\int\limits_{{x_1}}^x {d{x_2}\int\limits_{ - \frac{1}{2}}^1 {{{(x - {x_2})}^{ - \frac{1}{2}}}{{({x_2} - {x_1})}^{\frac{k}{2} - \frac{3}{2}}}dy_2.} }
\end{eqnarray*}

Evaluating the previous integrals we have that
\begin{eqnarray*}
  \left| {{K_k}(x,y;{x_1},{y_1})} \right| & \leq & {M^k}{(\sqrt \pi  )^{k - 1}}{\left(\frac{3}{2}\right)^{k - 1}}\frac{{{{(x - {x_1})}^{\frac{k}{2} - 1}}}}{{\Gamma (\frac{{k - 1}}{2})}}  \int\limits_0^1 {{\sigma ^{ - \frac{1}{2}}}{{(1 - \sigma )}^{\frac{k}{2} - \frac{3}{2}}}d\sigma} \\
  & =&{ (\sqrt \pi  } M{)^k}{\left(\frac{3}{2}\right)^{k - 1}}\frac{{{{(x - {x_1})}^{\frac{k}{2} - 1}}}}{{\Gamma (\frac{k}{2})}}, \hfill \\
\end{eqnarray*}
which proves Lemma 5. \qed

Now we are in a position to prove the final result of this paper, which gives us the Volterra property for the inverse of operator $\mathbb{L}$.

\textbf{Theorem 3.} The integral operator defined in the right hand of (25), i.e.
$$
{\mathbb{L}^{ - 1}}f(x,y) = \iint\limits_\Omega  {K(x,y;\,{x_1},{y_1})f({x_1},{y_1})d{x_1}d{y_1},}\eqno (33)
$$
has the Volterra property (it is almost continuous and quasi-nilpotent) in ${L_2}$($\Omega $).

\verb"Proof:"  Since the continuity of this operator follows from the fact that $K \in {L_2}(\Omega  \times \Omega )$. To prove this theorem, we only need to verify that operator ${\mathbb{L}^{ - 1}}$, defined by (33), is quasi-nilpotent, i.e.
$$
\mathop {\ell im}\limits_{n \to \infty } \left\| {{\mathbb{L}^{ - n}}} \right\|_0^{\frac{1}{n}} = 0,\eqno (34)
$$
where
$$
{\mathbb{L}^{ - n}} = {\mathbb{L}^{ - 1}}\left[ {{\mathbb{L}^{ - (n - 1)}}} \right]\,\,\,\,\,\,,\,\,\,\,n = 1,2,3, \ldots
$$

From (33), and by direct calculations, one can easily arrive at the following expression:
$$
{\mathbb{L}^{ - n}}f(x,y) = \iint\limits_\Omega  {{K_n}(x,y;\,{x_1},{y_1})f({x_1},{y_1})d{x_1}d{y_1}.}\eqno (35)
$$

Consequently, using the inequality of Schwarz and expression (32), from the representation (35) we obtain that
\begin{eqnarray*}
  \left\| {{\mathbb{L}^{ - n}}f} \right\|_0^2 &=& \iint\limits_\Omega  {{{\left| {{\mathbb{L}^{ - n}}f} \right|}^2}dxdy} = \iint\limits_\Omega  {{{\left[ {\iint\limits_\Omega  {{K_n}(x,y;{x_1},{y_1})f({x_1},{y_1})d{x_1}d{y_1}}} \right]}^2}dxdy  } \hfill \\
  & \leq & \iint\limits_\Omega  \left[ {\left( {\iint\limits_\Omega  {{{\left| {f\left( {{x_1},{y_1}} \right)} \right|}^2}d{x_1}d{y_1}}} \right)\,\left( {\iint\limits_\Omega  {{{\left| {{K_n}\left( {x,y;\,{x_1},{y_1}} \right)} \right|}^2}d{x_1}d{y_1}}} \right)} \right]dxdy \\
 & \leq & \left( {\frac{3}{2}\sqrt \pi  M} \right)^{2n}\frac{1}{{n\left( {n - 1} \right){\Gamma ^2}\left( {\frac{n}{2}} \right)}}\left\| f \right\|_0^2.
\end{eqnarray*}

From here we get
$$
{\left\| {{\mathbb{L}^{ - n}}} \right\|_0} \leq \left( {\frac{3}{2}\sqrt \pi  M} \right)^{n}\frac{1}{{\Gamma (1 + \frac{n}{2})}}.
$$
From the last equality one can state the validity of the equality (34) and Theorem 3 is proved.\qed

\textbf{Consequence 1.}  Problem B has the Volterra property.

\textbf{Consequence 2.} For any complex number $\lambda$, the equation
$$
\mathbb{L}u - \lambda u = f\eqno (36)
$$
is uniquely solvable for all $f \in {L_2}(\Omega )$.

Due to the invertibility of the operator $\mathbb{L}$, the unique solvability of the Eq.(36) is equivalent to the uniqueness of solution of the equation
$$
u - \lambda {\mathbb{L}^{ - 1}}u = {\mathbb{L}^{ - 1}}f,
$$
which is a second kind Volterra equation. This proves Consequence 2 of Theorem 3.

\section*{Competing interests}

The authors declare that they have no competing interests.

\section*{Authors' contributions}

The four authors have participated into the obtained results. The collaboration of each one cannot be separated in different parts of the paper. All of them have made substantial contributions to the theoretical results. The four authors have been involved in drafting the manuscript and revising it critically for important intellectual content.  All authors have given final approval of the version to be published.

\section*{Acknowledgement}
This research was partially supported by Ministerio de Ciencia e Innovaci\'{o}n-SPAIN, and FEDER, project MTM2010-15314 and KazNPU Rector's grant for 2013 .

\end{bmcformat}
\end{document}